\RequirePackage{amsmath}

\documentclass{svproc}

\usepackage{url}

\usepackage{acro}
\DeclareAcronym{gmres}{short=GMRES,long=Generalized Minimum Residual}
\DeclareAcronym{mgs}{short=MGS, long=Modified Gram-Schmidt}
\DeclareAcronym{cgsr}{short=CGSR,long=Classical Gram-Schmidt with Reorthogonalization}

\usepackage{amssymb}
\usepackage{algorithm}
\usepackage{algpseudocode}
\usepackage{graphicx}
\usepackage{booktabs}

\newcommand{\mat}[1]{\tens{#1}}
\newcommand{\realvecs}[1]{\mathbb{R}^{#1}}
\newcommand{\realmats}[2]{\mathbb{R}^{#1\times#2}}
\renewcommand{\epsilon}{\varepsilon}
\usepackage[binary-units=true]{siunitx}
\newcommand{\suitesparse}[1]{\texttt{#1}}

\usepackage{todonotes}

\begin{document}
\mainmatter              
\title{Improving the Performance of the GMRES Method using Mixed-Precision Techniques}
\titlerunning{Mixed-Precision GMRES}  
%
\author{Neil Lindquist\inst{1}\textsuperscript{[0000-0001-9404-3121]} \and Piotr Luszczek\inst{1} \and Jack Dongarra\inst{1}$^,$\inst{2}$^,$\inst{3}\textsuperscript{[0000-0003-3247-1782]}}
\authorrunning{N. Lindquist, P. Luszczek, and J. Dongarra} 
%
\tocauthor{Neil Lindquist, Piotr Luszczek, and Jack Dongarra}
\institute{University of Tennessee, Knoxville TN, USA
\email{\{nlindqu1,luszczek,dongarra\}@icl.utk.edu}
\and
Oak Ridge National Laboratory, Oak Ridge TN, USA
\and
University of Manchester, Manchester, UK}

\maketitle              

\begin{abstract}
The GMRES method is used to solve sparse, non-symmetric systems of
linear equations arising from many scientific applications.
The solver performance within a single node is memory bound,
due to the low arithmetic intensity of its computational kernels.
To reduce the amount of data movement, and thus,
to improve performance, we investigated the effect of using a mix of
single and double precision while retaining double-precision accuracy.
Previous efforts have explored reduced precision in the preconditioner,
but the use of reduced precision in the solver itself has received
limited attention. We found that GMRES only needs double precision in
computing the residual and updating the approximate solution to
achieve double-precision accuracy, although it must restart after each
improvement of single-precision accuracy.  This finding holds for the
tested orthogonalization schemes: Modified Gram-Schmidt (MGS) and
Classical Gram-Schmidt with Re-orthogonalization (CGSR).
Furthermore, our mixed-precision GMRES, when restarted at
least once, performed 19\% and 24\% faster on average than
double-precision GMRES for MGS and CGSR, respectively.
Our implementation uses generic programming techniques to ease the
burden of coding implementations for different data types. Our use of
the Kokkos library allowed us to exploit parallelism and optimize
data management. Additionally, KokkosKernels was used when producing
performance results. In conclusion, using a mix of single and
double precision in GMRES can improve performance while retaining
double-precision accuracy.
\keywords{Krylov subspace methods, mixed precision, linear algebra, Kokkos}
\end{abstract}

\section{Introduction}
The \ac{gmres} method~\cite{saadGMRESGeneralizedMinimal1986} is used for solving sparse, non-symmetric systems of linear equations arising from many applications~\cite[p. 193]{saadIterativeMethodsSparse2003}.
It is an iterative, Krylov subspace method that constructs an orthogonal basis by Arnoldi's procedure~\cite{arnoldiPrincipleMinimizedIteration1951} then finds the solution vector in that subspace such that the resulting residual is minimized.
One important extension of \ac{gmres} is the introduction of restarting,
whereby, after some number of iterations, \ac{gmres} computes
the solution vector, then starts over with an empty Krylov subspace and
the newly computed solution vector as the new initial guess.
This limits the number of basis vectors required for the Krylov subspace thus
reducing storage and the computation needed to orthogonalize each new vector.
On a single node system, performance of \ac{gmres} is bound by the main
memory bandwidth due to the low arithmetic intensity of its
computational kernels.  We investigated the use of a mix of single and
double floating-point precision to reduce the amount of data that needs
to be moved across the cache hierarchy, and thus improve the
performance, while trying to retain the accuracy that may be achieved by
a double-precision implementation of \ac{gmres}.
We utilized the iterative nature of \ac{gmres}, particularly when
restarted, to overcome the increased round-off errors introduced by
reducing precision for some computations.

The use of mixed precision in solving linear systems has long been established in the form of iterative refinement for dense linear systems~\cite{wilkinsonRoundingErrorsAlgebraic1963}, which is an effective tool for increasing performance~\cite{buttariMixedPrecisionIterative2007}.
However, research to improve the performance of \ac{gmres} in this way has had limited scope.
One similar work implemented iterative refinement with single-precision Krylov solvers, including \ac{gmres}, to compute the error corrections~\cite{anztMixedPrecisionIterative2010}.
However, that work did not explore the configuration of \ac{gmres} and tested only a limited set of matrices.
Recent work by Gratton et al.\@ provides detailed theoretical results for mixed-precision \ac{gmres}~\cite{grattonExploitingVariablePrecision2020}; although, they focus on non-restarting \ac{gmres} and understanding the requirements on precision for each inner-iteration to converge as if done in uniform, high precision.
Another approach is to use reduced precision only for the preconditioner~\cite{giraudMixedprecisionPreconditionersParallel2008}.
One interesting variant of reduced-precision preconditioners is to use a single-precision \ac{gmres} to precondition a double-precision \ac{gmres}~\cite{baboulinAcceleratingScientificComputations2008}.

In this paper, we focus on restarted \ac{gmres} with left
preconditioning and with one of two orthogonalization schemes: \ac{mgs}
or \ac{cgsr}, as shown in Alg.~\ref{alg:gmres}.  The algorithm contains
the specifics of the \ac{gmres} formulation that we used.
\begin{algorithm}[tbp]
	\caption{Restarted GMRES with left preconditioning~\cite{saadIterativeMethodsSparse2003}}
	\label{alg:gmres}
	\begin{algorithmic}[1]
		\State \(\mat{A} \in \realmats{n}{n},
		\quad\vec{x_0}, \vec{b} \in \realvecs{n},
		\quad \mat{M}^{-1} \approx \mat{A}^{-1}\)
		\For{\(k = 1, 2, \dots\)}
			\State{\(\vec{z_k} \gets \vec{b} - \mat{A}\vec{x_k}\)}
			\Comment{compute residual}
			\State{If \(\|\vec{z_k}\|_2\) is small enough, stop}
			\label{alg:gmres-low-prec-start}
			\State{\(\vec{r_k} \gets \mat{M}^{-1}\vec{z_k}\)}

			\State{\(\beta\gets \|\vec{r_k}\|_2,
					\quad s_0 \gets \beta,
					\quad\vec{v_1}\gets \vec{r}_k/\beta,
					\quad\mat{V_1} \gets [\vec{v_1}]\)}
			\State \(j \gets 0\)

			\Loop{ until the restart condition is met}
				\State \(j \gets j + 1\)
				\State{\(\vec{w} \gets \mat{M}^{-1}\mat{A}\vec{v_j}\)}
				\State{\(\vec{w}, h_{1,j}, \dots, h_{j,j}
						\gets\mathrm{orthogonalize}(\vec{w}, \mat{V_j})\)}
				\Comment{\ac{mgs} or \ac{cgsr}}
				\State{\(h_{j+1,j} \gets \|\vec{w}\|_2\)}
				\State{\(\vec{v_{j+1}} \gets \vec{w}/h_{j+1,j}\)}
				\State\(\mat{V_{j+1}} \gets [\mat{V_j}, \vec{v_{j+1}}]\)

				\For{\(i = 1, \dots, j-1\)}
				\State \(
					\begin{bmatrix} h_{i,j} \\ h_{i+1,j} \end{bmatrix}
					\gets
					\begin{bmatrix} \alpha_i & \beta_i \\
					-\beta_i & \alpha_i \end{bmatrix}
					\times
					\begin{bmatrix} h_{i,j} \\ h_{i+1,j} \end{bmatrix}
				\) \Comment{apply Givens rotation}
				\EndFor
				\State \(
					\begin{bmatrix} \alpha_j \\ \beta_j \end{bmatrix}
					\leftarrow
					\text{\textsf{rotation\_matrix}}
					\left( \begin{bmatrix} h_{j,j} \\ h_{j+1,j} \end{bmatrix} \right)
				\) \Comment{form $j$-th Givens rotation}

				\State \(
					\begin{bmatrix} s_{j} \\ s_{j+1} \end{bmatrix}
					\leftarrow
					\begin{bmatrix} \alpha_j & \beta_j \\
					-\beta_j & \alpha_j \end{bmatrix}
					\times
					\begin{bmatrix} s_{j} \\ 0 	\end{bmatrix}
				\)
				\State \(
					\begin{bmatrix} h_{j,j} \\ h_{j+1,j} \end{bmatrix}
					\leftarrow
					\begin{bmatrix} \alpha_j & \beta_j \\
					-\beta_j & \alpha_j \end{bmatrix}
					\times
					\begin{bmatrix} h_{j,j} \\ h_{j+1,j} \end{bmatrix}
				\)
			\EndLoop

			\State{\(\mat{H} \gets \{h_{i,\ell}\}_{1\leq i,\ell \leq j},\quad
				\vec{s} \gets [s_1, \dots s_j]^T\)}
			\State \(\vec{u_k} \gets \mat{V_j}\mat{H}^{-1}\vec{s}\) \Comment{compute correction}
			\label{alg:gmres-low-prec-end}
			\State \(\vec{x_{k+1}} \gets \vec{x_k} + \vec{u_k}\) \Comment{apply correction}

		\EndFor
		\Statex
		\Procedure{\ac{mgs}}{\(\vec{w}, \mat{V_j}\)\relax}
			\State\([\vec{v_1}, \dots, \vec{v_{j}}] \gets \mat{V_{j}}\)
			\For{i = 1, 2, \dots, j}
				\State{\(h_{i,j} \gets \vec{w}\cdot \vec{v_i}\)}
				\State{\(\vec{w} \gets \vec{w} - h_{i,j}\vec{v_i}\)}
			\EndFor
			\State\Return \(\vec{w}, h_{1,j}, \dots, h_{j,j}\)
		\EndProcedure
		\Statex
		\Procedure{\ac{cgsr}}{\(\vec{w}, \mat{V_j}\)\relax}
			\State \(\vec{h} \gets \mat{V_j}^T\vec{w} \)
			\State \( \vec{w} \gets \vec{w}-\mat{V_j} \vec{h} \)
			\State \([h_{0,j}, \dots, h_{j,j}]^T \gets \vec{h}\)
			\State\Return \(\vec{w}, h_{1,j}, \dots, h_{j,j}\)
		\EndProcedure
	\end{algorithmic}
\end{algorithm}
\ac{mgs} is the usual choice for orthogonalization in \ac{gmres} due to its lower computational cost compared to other schemes~\cite{paigeResidualBackwardError2001}.
\ac{cgsr} is used less often in practice but differs in interesting ways from \ac{mgs}.
First, it retains good orthogonality relative to round-off
error~\cite{giraudLossOrthogonalityGramSchmidt2005}, which raises the
question of whether this improved orthogonality can be used to circumvent
some loss of precision.
Second, it can be implemented as matrix-vector multiplies, instead of
a series of dot-products used by \ac{mgs}.
Consequently, \ac{cgsr} requires fewer global reductions and may be a
better candidate when considering expanding the work to a distributed
memory setting.  Restarting is used to limit the storage and computation
requirements of the Krylov basis generated by
\ac{gmres}~\cite{Baker2003gmresimprove,Baker2005accrstrtgmres}.

\section{Numerics of Mixed Precision \ac{gmres}}
\label{sec:numerics}

To use mixed precision for improving the performance of \ac{gmres}, it
is important to understand how the precision of different parts of the
solver affects the final achievable accuracy.
First, for the system of linear equations \(\mat{A} \vec{x}=\vec{b}\);
\(\mat{A}\), \(\vec{b}\), and \(\vec{x}\) all must be stored in full precision because changes to these values change the problem being solved and directly affect the backward and forward error bounds.
Next, note that restarted \ac{gmres} is equivalent to iterative refinement where the error correction is computed by non-restarted \ac{gmres}.
Hence, adding the error correction to the current solution must be done
in full precision to prevent \(\vec{x}\) from suffering round-off to reduced precision.
Additionally, full precision is critical for the computation of residual
\(\vec{r}=\mat{A}\vec{x}-\vec{b}\) because it is used to compute the error correction and
is computed by subtracting quantities of similar value.
Were the residual computed in reduced precision, the maximum error that could be corrected is limited by the accuracy used for computing the residual vector~\cite{carsonNewAnalysisIterative2017}.

Next, consider the effects of reducing precision in the computation of the error correction.
Note that for stationary iterative refinement algorithms, it has long been known
that reduced precision can be used in this way while still achieving
full accuracy~\cite{wilkinsonRoundingErrorsAlgebraic1963}, which, to
some extent, can be carried over to the non-stationary correction of \ac{gmres}.
The converge property derives from the fact that if each restart
\({i = 1, \dots, k}\) computes an update, \(\vec{u}_i\), fulfilling
\({\|\vec{r}_i\|} = {\|\vec{r}_{i-1} - \mat{A}\vec{u}_i\|} \leq
{\delta\|\vec{r}_{i-1}\|}\) for some error reduction
\(\delta<1\), then after $k$ steps we get \({\|\vec{r}_k\|}\leq
{\delta^k\|\vec{r}_0\|}\)~\cite{anztMixedPrecisionIterative2010}.
Thus, reducing the accuracy of the error-correction to single precision
does not limit the maximum achievable accuracy.
Furthermore, under certain restrictions on the round-off errors of
the performed operations, non-restarted \ac{gmres} behaves as if the
arithmetic was done
exactly~\cite{grattonExploitingVariablePrecision2020}.
Therefore, when restarted frequently enough, we hypothesize that
mixed-precision \ac{gmres} should behave like the double-precision
implementation.

\section{Restart Strategies}
\label{sec:restart}

Restart strategies are important to the convergence.
In cases when limitations of the memory use require a \ac{gmres} restart
before the accuracy in working precision is reached, the restart
strategy needs no further consideration.
However, if mixed-precision \ac{gmres} may reach the reduced precision's accuracy before the iteration limit, it is important to have a strategy to restart early.
But restarting too often will reduce the rate of convergence because
improvement is related to the Arnoldi process's approximation of the eigenvalues
of \(\mat{A}\), which are discarded when \ac{gmres} restarts~\cite{vandervorstSuperlinearConvergenceBehaviour1993}.
As a countermeasure, we propose four possible approaches for robust convergence monitoring and restart initiation.

There are two points to note before discussing specific restart strategies.
First, the choice of orthogonalization scheme is important to consider,
because some Krylov basis vectors usually become linearly dependent when
\ac{gmres} reaches the working precision accuracy, e.g.,
\ac{mgs}~\cite{paigeResidualBackwardError2001}, while other methods
remain nearly orthogonal, e.g., \ac{cgsr}~\cite{giraudLossOrthogonalityGramSchmidt2005}.
Second, the norm of the Arnoldi residual, the residual for
\ac{gmres}'s least-squares problem, approximates the norm of the residual
of the original preconditioned linear system of equations and is
computed every iteration when using Givens rotations to solve the
least-squares
problem (\(s_{j+1}\) in Alg.~\ref{alg:gmres})~\cite[Proposition~6.9]{saadIterativeMethodsSparse2003}.
However, this approximation only monitors the least-squares problem and is not guaranteed to be accurate after reaching working precision~\cite{greenbaumNumericalBehaviourModified1997}.
The explanation is unknown, but it has been noted that the Arnoldi
residual typically decreases past the true residual if and only if
independent vectors continue to be added to the Krylov basis.
Hence, the choice of orthogonalization scheme must be considered
when using restarts based on the Arnoldi residual norm.

Our first restart strategy derives from the observation that the number
of iterations before the convergence stalls appears to be roughly
constant after each restart.
See Sec.~\ref{sec:results-convergence} for numerical examples.
While this does not alleviate the issue determining the appropriate point for the first restart, this can be used for subsequent restarts either to trigger the restart directly or as a heuristic for when to start monitoring other, possibly expensive, metrics.

The second restart strategy is to monitor the approximate preconditioned
residual norm until it drops below a given threshold, commonly related to
the value after the prior restart.
The simplest threshold is a fixed, scalar value.
Note that if the approximated norm stops decreasing, such as for \ac{mgs}, this criterion will not be met until \ac{gmres} is restarted.
Thus, the scalar thresholds must be carefully chosen when using \ac{mgs}.
More advanced threshold selection may be effective, but we have not
explored any yet.

Inspired by the problematic case of the second strategy,
the third strategy is to detect when the Arnoldi residual norm stops improving.
Obviously, this approach is only valid if the norm stops decreasing when \ac{gmres} has stalled.
Additionally, \ac{gmres} can stagnate during normal operation, resulting in iterations of little or no improvement, which may cause premature restarts.

The final strategy is to detect when the orthogonalized basis becomes linearly dependent.
This relates to the third strategy but uses a different approach.
For the basis matrix \(\mat{V}_k\) computed in the \(k\)th inner iteration, let \(\mat{S}_k=(\mat{I}+\mat{U}_k)^{-1}\mat{U}_k\), where \(\mat{U}_k\) is the strictly upper part of \(\mat{V}_k^H\mat{V}_k\)~\cite{paigeUsefulFormUnitary2009}.
Then, the basis is linearly dependent if and only if \(\|\mat{S}_k\|_2=1\).
It has been conjectured that \acs{mgs}-\acs{gmres} converges to machine precision when the Krylov basis loses linear independence~\cite{paigeResidualBackwardError2001,paigeModifiedGramSchmidtMGS2006}.
This matrix can be computed incrementally, appending one column per inner iteration, requiring \(2nk + 2k^2\) FLOP per iteration.
Estimating the 2-norm for a \ac{gmres}-iteration with \(i\) iterations
of the power method requires an additional \(i(2k^2+3k)\) FLOP by
utilizing the strictly upper structure of the matrix.

\section{Experimental Results}

First, Sec.~\ref{sec:results-convergence} shows accuracy and rate of convergence results to verify the results in Sec.~\ref{sec:numerics} and to better understand the strategies proposed in Sec.~\ref{sec:restart}.
Next, Sec.~\ref{sec:results-performance} compares the performance of our mixed-precision approach and double-precision \ac{gmres}.
Based on the results in Sec.~\ref{sec:numerics}, we focused on computing the residual and updating \(\vec{x}\) in double-precision and computing everything else in single-precision.
This choice of precisions has the advantage that it can be implemented
using uniform-precision kernels and only casting the residual to
single-precision and the error-correction back to double precision.
Note that in this approach, the matrix is stored twice, once in single-precision and once in double-precision; this storage requirement may be able to be improved by storing the high-order and low-order bytes of the  double-precision matrix in separate arrays~\cite{grutzmacherModularPrecisionFormat2018}.

Matrices were stored in Compressed Sparse Row format and preconditioned with incomplete LU without fill in.
So, the baseline, double-precision solver requires \(24n_{nz} + 8nm + 28n + 8m^2 + O(m)\) bytes while the mixed-precision solver requires \(24n_{nz} + 4nm + 32n + 4m^2 + O(m)\) where \(n_{nz}\) is the number of matrix nonzero elements, \(n\) is the number of matrix rows, and \(m\) is the maximum number of inner iterations per restart.
All of the tested matrices came from the SuiteSparse
collection~\cite{davisUniversityFloridaSparse2011} and entries of the
solution vectors were independently drawn from a uniform distribution between 0 and 1.

We used two implementations of \ac{gmres}: a configurable one for exploring the effect various factors have on the rate of convergence, and an optimized one for testing a limited set of factors for performance.
Both implementations are based on version 2.9.00 of the Kokkos performance portability library~\cite{edwardsKokkosEnablingManycore2014}.
The OpenMP backend was used for all tests.
Furthermore, for performance results we used the KokkosKernels library, with Intel's MKL where supported, to ensure that improvements are compared against a state-of-the-art baseline.
The rate of convergence tests were implemented using a set of custom, mixed-precision kernels for ease of experimentation.

All experiments were run on a single node with two sockets, each
containing a ten-core Haswell processor, for a total of twenty-cores and
\SI{25}{\mebi\byte} of combined Level 3 cache.
Performance tests were run with Intel C\texttt{++} Compiler version 2018.1, Intel MKL version 2019.3.199, and Intel Parallel Studio Cluster Edition version 2019.3.
The environment variables controlling OpenMP were set to:
\verb|OMP_NUM_THREADS=20|, \verb|OMP_PROC_BIND=spread|, and
\verb|OMP_PROC_BIND=places|.

\subsection{Measurement of the Rate of Convergence}
\label{sec:results-convergence}

\begin{figure}[tbp]
	\centering
	\includegraphics{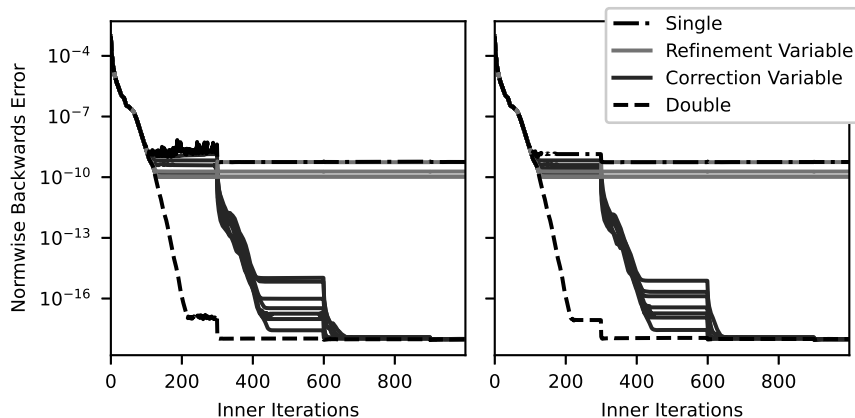}
	\caption{Rate of convergence results for the \suitesparse{airfoil\_2d} matrix when restarting every 300 iterations for \ac{mgs} (left) and \ac{cgsr} (right) orthogonalization schemes}
	\label{fig:convergence-pervariable-airfoil_2d}
\end{figure}

To verify the analysis of Sec.~\ref{sec:numerics}, we first
demonstrate that each variable behaves as predicted when stored in
single-precision, while the rest of the solver components remain in double-precision.
Figure~\ref{fig:convergence-pervariable-airfoil_2d} shows the normwise backward error after each inner iteration as if the solver had terminated, for \ac{gmres} solving a linear system for the \suitesparse{airfoil\_2d} matrix.
This matrix has \num{14214} rows, \num{259688} nonzeros, and a condition of \num{1.8e6}.
In the figure, the ``Refinement Variables'' include the matrix when used for the residual, the right-hand side, the solution, and the vector used to compute the non-preconditioned residual; the ``Correction Variables'' include the matrix when used to compute the next Krylov vector, the non-preconditioned residual, the Krylov vector being orthogonalized, the orthogonal basis, the upper triangular matrix from the orthogonalization process, and the vectors to solve the least-squares problems with Givens rotations.
The convergence when storing the preconditioner in single-precision was
visually indistinguishable from the double-precision baseline and
omitted from the figure for the sake of clarity.
Each solver was restarted after 300 iterations.
All of the solvers behaved very similarly until single-precision accuracy was reached, where all of the solvers, except double-precision, stopped improving.
After restarting, the solvers with reduced precision inside the error correction started improving again and eventually reached double-precision accuracy; however, the solvers with reduced precision in computing the residual or applying the error correction were unable to improve past single-precision accuracy.

\begin{figure}[tbp]
	\centering
	\includegraphics{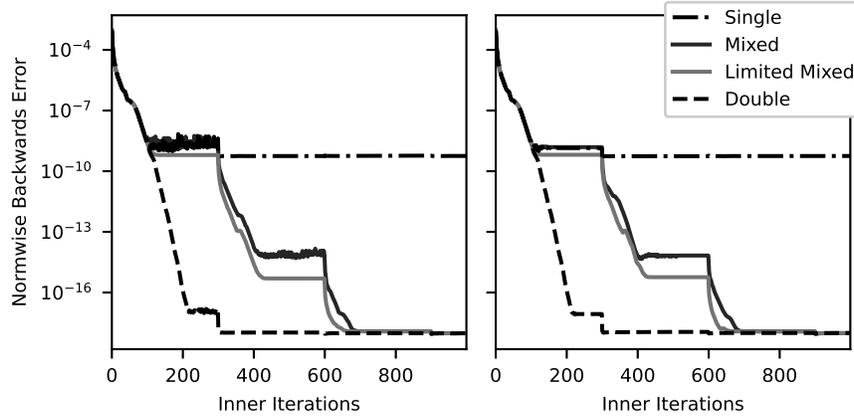}
	\caption{Rate of convergence results for the \suitesparse{airfoil\_2d} matrix when restarting every 300 iterations for \ac{mgs} (left) and \ac{cgsr} (right) orthogonalization schemes}
	\label{fig:convergence-airfoil_2d}
\end{figure}

The convergence test was repeated with two mixed-precision solvers that use reduced precision for multiple variables.
The first used double precision only for computing the residual and error correction, i.e., using single precision for lines~\ref{alg:gmres-low-prec-start}-\ref{alg:gmres-low-prec-end} of Alg.~\ref{alg:gmres}.
The second was more limited, using single precision only to store \(\mat{A}\) for computing the next Krylov vector, the preconditioner \(\mat{M}^{-1}\), and the Krylov basis \(\mat{V_j}\) from Alg.~\ref{alg:gmres}; these three variables make up most of the data that can be stored in reduced precision.
Figure~\ref{fig:convergence-airfoil_2d} shows the normwise backward error after each inner iteration for single, double, and mixed precisions solving a linear system for the \suitesparse{airfoil\_2d} matrix.
After restarting, both mixed-precision \ac{gmres} implementations were able to resume improvement and achieve double-precision accuracy.
This ability to converge while using reduced precision occurred for all of the matrices tested, as can be seen in Sec.~\ref{sec:results-performance}.
Note that while limiting the use of mixed precision can increase the amount of improvement achieved before stalling, this improvement is limited and does not reduce the importance of appropriately restarting.
Additionally, the limited mixed-precision implementation requires several mixed-precision kernels, while the fully mixed-precision implementation can be implemented using uniform-precision kernels by copying the residual to single-precision and copying the error-correction back to double-precision.

One interesting observation was that the number of iterations before improvement stalled was approximately the same after each restart.
Table~\ref{tab:restarts-stall-point} displays the number of iterations before stalling after the first three restarts in the mixed-precision \ac{mgs}-\ac{gmres}.
\begin{table}[tbp]
	\caption{Number of iterations before the improvement stalls in mixed-precision \ac{mgs}-\ac{gmres}}
	\label{tab:restarts-stall-point}
	\centering
	\begin{tabular}{r r r r r}
		       & {Iterations}  & {Iterations}    & {Iterations}    & {Iterations} \\
		Matrix & {per Restart} & {for 1st Stall} & {for 2nd Stall} & {for 3rd Stall} \\
		\midrule
		\suitesparse{airfoil\_2d} & 300 & 137 & 141 & 142 \\
		\suitesparse{big} & 500 & 360 & 352 & 360 \\
		\suitesparse{cage11} & 20 & 7 & 7 & 8 \\
		\suitesparse{Goodwin\_040} & 1250 & 929 & 951 & 924 \\
		\suitesparse{language} & 75 & 23 & 21 & 21 \\
		\suitesparse{torso2} & 50 & 28 & 27 & 25 \\
	\end{tabular}
\end{table}
Stalling was defined here to be the Arnoldi residual norm improving by less than a factor of \num{1.001} on the subsequent 5\% of inner iterations per restart.
This behavior appears to hold for \ac{cgsr} too but was not quantified
because stalled improvement cannot be detected in the Arnoldi residual for \ac{cgsr}.

\begin{figure}[tbp]
	\centering
	\includegraphics{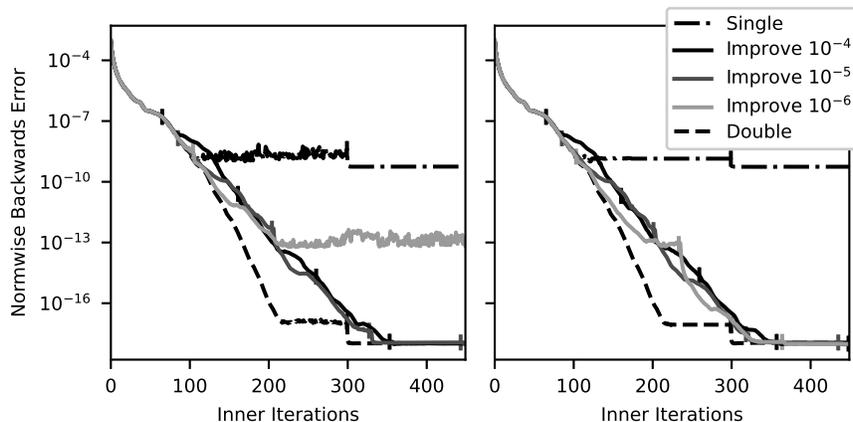}
	\caption{Rate of convergence results for the \suitesparse{airfoil\_2d} matrix when restarting mixed-precision \ac{gmres} after a fixed improvement in the Arnoldi residual norm  for \ac{mgs} (left) and \ac{cgsr} (right) orthogonalization schemes, with vertical ticks to indicate when restarts occurred}
	\label{fig:restart-fixedthresh-airfoil_2d}
\end{figure}

\begin{figure}[tbp]
	\centering
	\includegraphics{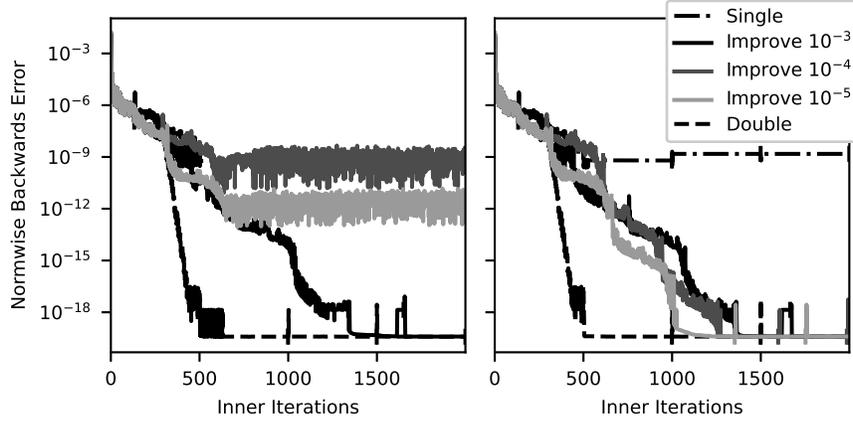}
	\caption{Rate of convergence results for the \suitesparse{big} matrix when restarting mixed-precision \ac{gmres} after a fixed improvement in the Arnoldi residual norm  for \ac{mgs} (left) and \ac{cgsr} (right) orthogonalization schemes, with vertical ticks to indicate when restarts occurred}
	\label{fig:restart-fixedthresh-big}
\end{figure}

\begin{figure}[tbp]
	\centering
	\includegraphics{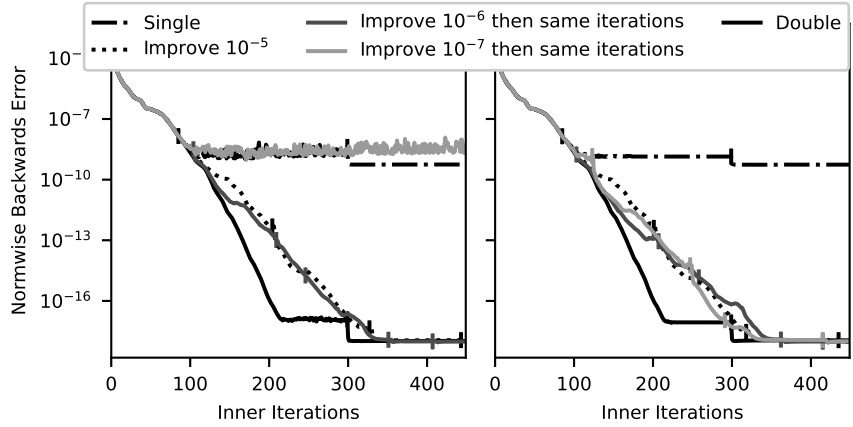}
	\caption{Rate of convergence results for the \suitesparse{airfoil\_2d} matrix when restarting mixed-precision \ac{gmres} after a fixed improvement in the Arnoldi residual norm for the first iteration and the same number of iterations thereafter for \ac{mgs} (left) and \ac{cgsr} (right) orthogonalization schemes, with vertical ticks to indicate when restarts occurred.  The rate of convergence using just a fixed improvement threshold of \(10^{-5}\) is added for comparison's sake}
	\label{fig:restart-fixedthresh+repeatiter-airfoil_2d}
\end{figure}

\begin{figure}[tbp]
	\centering
	\includegraphics{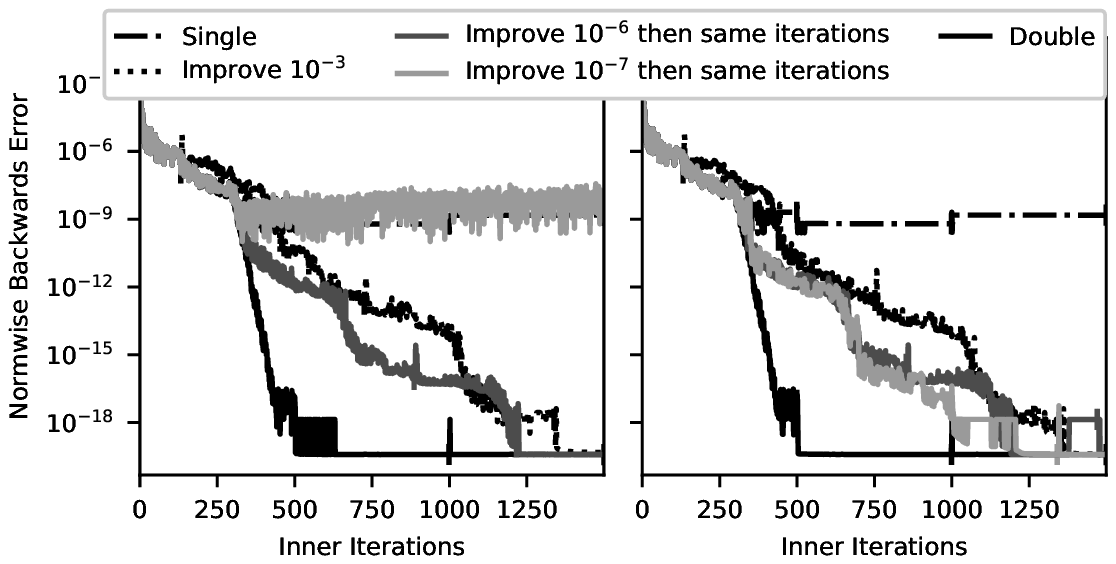}
	\caption{Rate of convergence results for the \suitesparse{big} matrix when restarting mixed-precision \ac{gmres} after a fixed improvement in the Arnoldi residual norm for the first iteration and the same number of iterations thereafter for \ac{mgs} (left) and \ac{cgsr} (right) orthogonalization schemes, with vertical ticks to indicate when restarts occurred.  The rate of convergence using just a fixed improvement threshold of \(10^{-5}\) is added for comparison's sake}
	\label{fig:restart-fixedthresh+repeatiter-big}
\end{figure}

Next, restart strategies based on the Arnoldi residual norm were tested.
First, Fig.~\ref{fig:restart-fixedthresh-airfoil_2d} shows the convergence when restarted after a fixed improvement.
Note that for \ac{mgs}, when the threshold is too ambitious, mixed-precision \ac{gmres} will stall because of roundoff error before reaching the threshold, at which point the approximated norm stops decreasing.
However, the choice of restart threshold becomes problematic when considering multiple matrices.
Figure~\ref{fig:restart-fixedthresh-big} shows the same test applied to the \suitesparse{big} matrix, which has \num{13209} rows, \num{91465} nonzeros, and an L2 norm of \num{4.4e7}.
Note that the successful threshold with the most improvement per restart is two orders of magnitude less improvement per restart than that of \suitesparse{airfoil\_2d}.
Next, Fig.~\ref{fig:restart-fixedthresh+repeatiter-airfoil_2d} uses the first restart's iteration count as the iteration limit for the subsequent restarts when solving the \suitesparse{airfoil\_2d} system.
Because only the choice of the first restart is important, a more ambitious threshold was chosen than for Fig.~\ref{fig:restart-fixedthresh-airfoil_2d}.
Note that, except for when the first restart was not triggered, this two-staged approach generally performed a bit better than the simple threshold.
Figure~\ref{fig:restart-fixedthresh+repeatiter-big} shows the mixed restart strategy for the \suitesparse{big} matrix.
Note how the same thresholds were used for the \suitesparse{big} test as the \suitesparse{airfoil\_2d} test but were still able to converge and outperform the matrix-specific, scalar threshold.
This two-part strategy appears to behave more consistently than the simple threshold.

\begin{figure}[tbp]
	\centering
	\includegraphics{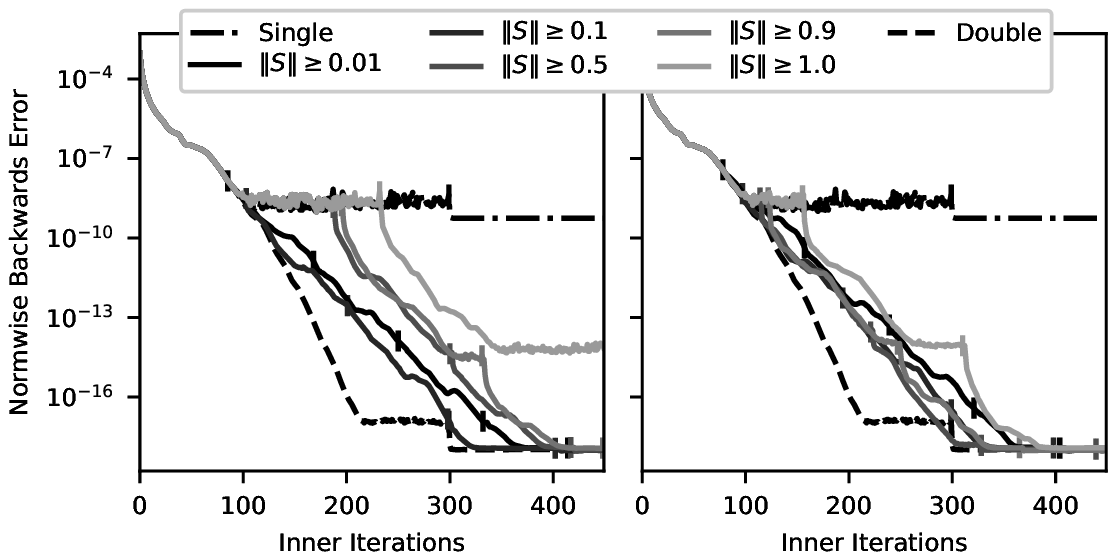}
	\caption{Rate of convergence results for the \suitesparse{airfoil\_2d} matrix when restarting mixed-precision \ac{gmres} based on the spectral norm (left) or Frobenius norm (right) of the \(S\) matrix, for \ac{mgs} orthogonalization, with vertical ticks to indicate when restarts occurred}
	\label{fig:restart-orthloss}
\end{figure}

Finally, we tested restarts based on the loss of orthogonality in the basis.
Because \ac{cgsr} retains a high degree of orthogonality, this strategy was only tested with \ac{mgs}-\ac{gmres}.
Fig.~\ref{fig:restart-orthloss} shows the rate of convergence when restarting based on the norm of the \(S\) matrix.
The spectral norm was computed using 10 iterations of the power method.
Additionally, the Frobenius norm was tested as a cheaper alternative to the spectral norm, although it does not provide the same theoretical guarantees.
Interestingly, when using the spectral norm, a norm of even \num{0.5} was not detected until improvement had stalled for a noticeable period.
Note that even the Frobenius norm, which is an upper bound on the spectral norm, did not reach \num{1} until, visually, improvement had stalled for a few dozen iterations.
The cause of this deviation from the theoretical results~\cite{paigeUsefulFormUnitary2009} is unknown.

\subsection{Performance}
\label{sec:results-performance}

Finally, we looked at the effect of reduced precision on performance.
Additionally, in testing a variety of matrices, these tests provide further support for some of the conclusions from Sec.~\ref{sec:results-convergence}.
The runtimes include the time spent constructing the preconditioner and making any copies of the matrix.
In addition to comparing the performance of mixed- and double-precision \ac{gmres}, we tested the effect of reducing the precision of just the ILU preconditioner.

\begin{figure}[tbp]
	\centering
	\includegraphics{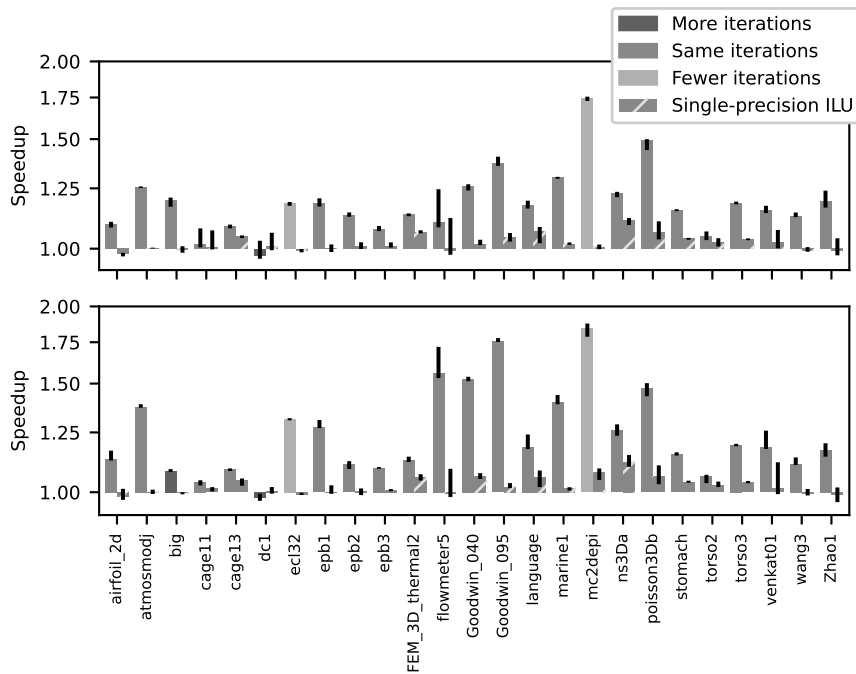}
	\caption{Speedup of the median runtime out of five tests for mixed-precision versus double-precision restarted in half the number of iterations needed for double-precision, for \ac{mgs} (top) and \ac{cgsr} (bottom) orthogonalization schemes, with error bars indicating the minimum and maximum speedups}
	\label{fig:performance-forced-restart}
\end{figure}

We first tested the performance improvement when other constraints force \ac{gmres} to restart more often than required by mixed precision.
For each of the tested systems, we computed the number of iterations for the double-precision solver to reach a backward error of \(10^{-10}\).
Then, we measured the runtime for each solver to reach a backward error
of \(10^{-10}\) when restarting after half as many iterations.
All but 3 of the systems took the same number of iterations for \ac{mgs}; two systems took fewer iterations for mixed precision (\suitesparse{ecl32} and \suitesparse{mc2depi}), while one system took more iterations for mixed precision (\suitesparse{dc1}).
\ac{cgsr} added one additional system that took more iterations for mixed precision (\suitesparse{big}).
Figure~\ref{fig:performance-forced-restart} shows the speedup of the mixed-precision implementation and the single-precision ILU implementation relative to the baseline implementation for each of the tested matrices.
For the mixed-precision implementation, the geometric mean of the speedup was 19\% and 24\% for \ac{mgs} and \ac{cgsr}, respectively.
For the single-precision ILU implementation, those means were both 2\%.

\begin{figure}[tbp]
	\centering
	\includegraphics{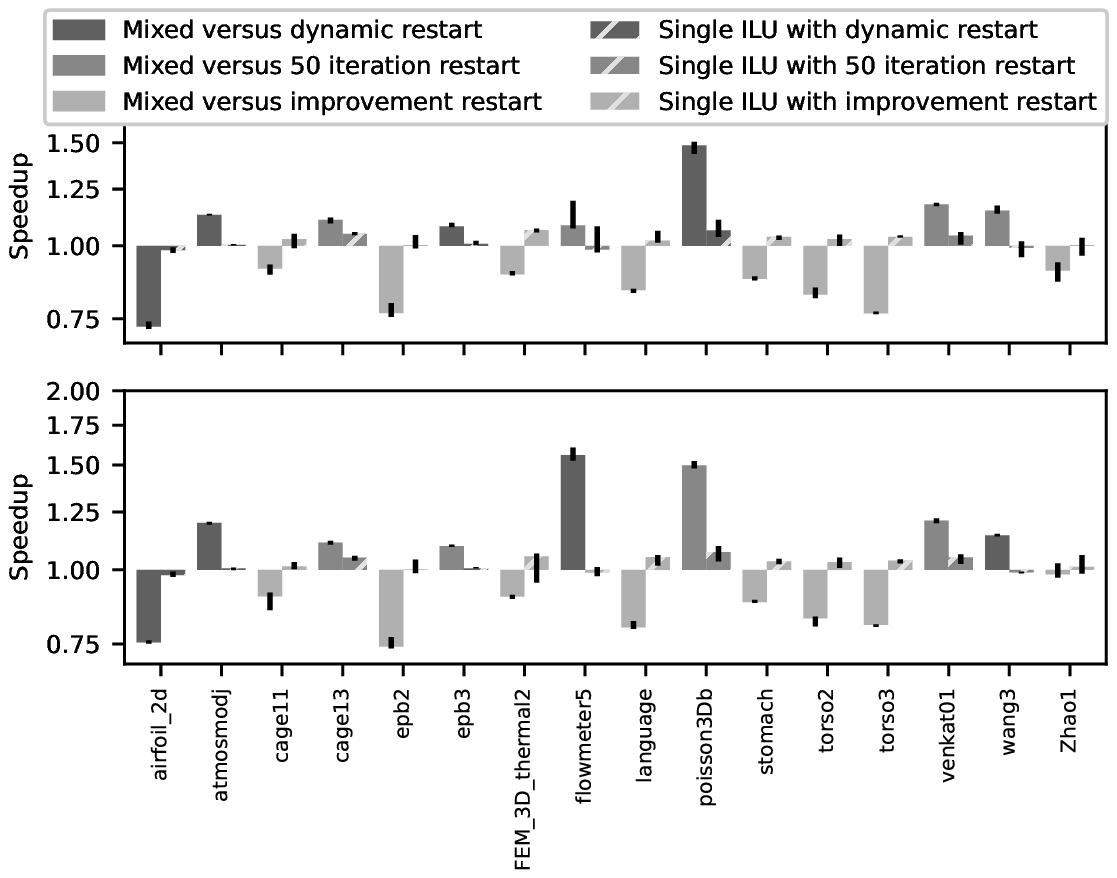}
	\caption{Speedup of the median runtime out of five tests for mixed-precision versus double-precision restarted after 50 iterations or an improvement in the Arnoldi residual, for \ac{mgs} (top) and \ac{cgsr} (bottom) orthogonalization schemes, with error bars indicating the minimum and maximum speedups}
	\label{fig:performance-restart-50}
\end{figure}

The second set of performance tests show what happens when \ac{gmres} is not forced to restart often enough for mixed precision.
All of the matrices from Fig.~\ref{fig:performance-forced-restart} that were restarted after fewer than 50 iterations were tested again, except they were restarted after 50 iterations.
For mixed-precision \ac{gmres}, the first restart could additionally be triggered by an improvement in the Arnoldi residual by a factor of \(10^{-6}\) and subsequent restarts were triggered by reaching the number of inner-iterations that caused the first restart.
To ensure the mixed-precision solver was not given any undue advantage,
the other two solvers' performance was taken as the best time from three
restart strategies: (1) the same improvement-based restart trigger as
mixed-precision \ac{gmres}; (2) after 50 iterations, or (3)  after an
improvement in the Arnoldi residual by a factor of \(10^{-8}\).
Figure~\ref{fig:performance-restart-50} shows the new performance results.
For the mixed-precision implementation, the geometric mean of the speedup was
-4\% and 0\% for \ac{mgs} and \ac{cgsr}, respectively.
For the single-precision ILU implementation, those means were 2\% and 1\% respectively.
The matrices for which the mixed-precision implementation performed
worse were exactly the matrices that did not require restarting when solved by the double-precision implementation.

\section{Conclusion}
As a widely used method for solving sparse, non-symmetric systems of linear equations, it is important to explore ways to improve the performance of \ac{gmres}.
Towards this end, we experimented with the use of mixed-precision techniques to reduce the amount of data moved across the cache hierarchy to improve performance.
By viewing \ac{gmres} as a variant of iterative refinement,
we found that \ac{gmres} was still able to achieve the accuracy of a
double-precision solver while using our proposed techniques of
mixed-precision and restart initiation.
Furthermore, we found that the algorithm, with our proposed
modifications, delivers improved performance when the
baseline implementation already requires restarting for all but one problem,
even compared to storing the preconditioner in single precision.
However, our approach reduced performance when the baseline
did not require restarting, at least for problems that require less
than 50 inner iterations.

There are a few directions in which this work can be further extended.
The first direction is to expand the implementation to a variety of
systems that are different from a single CPU-only node.
For example, GPU accelerators provide a significantly higher performance benefit than
CPUs but involve a different trade-off between computational units, memory
hierarchy, and kernel launch overheads.
Thus, it would be beneficial to explore the use of mixed precision in
\ac{gmres} on these systems.
Also important are the distributed memory, multi-node
systems that are used to solve problems too large to be computed
efficiently on a single node.
In these solvers, the movement of data across the memory hierarchy
becomes less important because of the additional cost of moving data between nodes.
A related direction is to explore the use of mixed-precision techniques to improve variants of \ac{gmres}.
One particularly important class of variants is communication-avoiding
and pipelined renditions for distributed systems, which use alternative formulations to reduce the amount of inter-node communication.
The last major direction is to explore alternative techniques to reduce data movement.
This can take many forms, including alternative floating-point
representations, such as half-precision, quantization, or Posits~\cite{gustafsonBeatingFloatingPoint2017};
alternative data organization, such as splitting the high- and low-order
bytes of double-precision~\cite{grutzmacherModularPrecisionFormat2018};
or applying data compression, such as SZ~\cite{liangErrorcontrolledLossyCompression2018} or
ZFP~\cite{lindstromFixedrateCompressedFloatingpoint2014}.

\subsubsection*{Acknowledgments}

This material is based upon work supported by the University of
Tennessee grant MSE E01-1315-038 as Interdisciplinary Seed funding and
in part by UT Battelle subaward 4000123266.
This material is also based upon work supported by the National Science
Foundation under Grant No.~2004541.

%
%
\bibliographystyle{spmpsci}
\bibliography{bibliography}

\begin{thebibliography}{10}
\providecommand{\url}[1]{{#1}}
\providecommand{\urlprefix}{URL }
\expandafter\ifx\csname urlstyle\endcsname\relax
  \providecommand{\doi}[1]{DOI~\discretionary{}{}{}#1}\else
  \providecommand{\doi}{DOI~\discretionary{}{}{}\begingroup
  \urlstyle{rm}\Url}\fi

\bibitem{anztMixedPrecisionIterative2010}
Anzt, H., Heuveline, V., Rocker, B.: Mixed precision iterative refinement
  methods for linear systems: Convergence analysis based on {Krylov} subspace
  methods.
\newblock In: Proceedings of the 10th International Conference on Applied
  Parallel and Scientific Computing - Volume 2, {{PARA}}'10, pp. 237--247.
  {Springer-Verlag}, {Berlin, Heidelberg} (2010).
\newblock \doi{10.1007/978-3-642-28145-7_24}

\bibitem{arnoldiPrincipleMinimizedIteration1951}
Arnoldi, W.E.: The principle of minimized iteration in the solution of the
  matrix eigenvalue problem.
\newblock Quart. Appl. Math. \textbf{9}, 17--29 (1951).
\newblock \doi{10.1090/qam/42792}

\bibitem{baboulinAcceleratingScientificComputations2008}
Baboulin, M., Buttari, A., Dongarra, J., Kurzak, J., Langou, J., Langou, J.,
  Luszczek, P., Tomov, S.: Accelerating scientific computations with mixed
  precision algorithms.
\newblock CoRR \textbf{abs/0808.2794} (2008).
\newblock \doi{10.1016/j.cpc.2008.11.005}

\bibitem{Baker2003gmresimprove}
Baker, A.H.: On improving the performance of the linear solver restarted
  {GMRES}.
\newblock Ph.D. thesis, University of Colorado (2003)

\bibitem{Baker2005accrstrtgmres}
Baker, A.H., Jessup, E.R., Manteuffel, T.: A technique for accelerating the
  convergence of restarted {GMRES}.
\newblock SIAM J. Matrix Anal. Appl. \textbf{26}(962) (2005).
\newblock \doi{10.1137/S0895479803422014}

\bibitem{buttariMixedPrecisionIterative2007}
Buttari, A., Dongarra, J., Langou, J., Langou, J., Luszczek, P., Kurzak, J.:
  Mixed precision iterative refinement techniques for the solution of dense
  linear systems.
\newblock Int. J. High Perform. Comput. Appl. \textbf{21}(4), 457--466 (2007).
\newblock \doi{10.1177/1094342007084026}

\bibitem{carsonNewAnalysisIterative2017}
Carson, E., Higham, N.J.: A new analysis of iterative refinement and its
  application to accurate solution of ill-conditioned sparse linear systems.
\newblock SIAM J. Sci. Comput. \textbf{39}(6), A2834--A2856 (2017).
\newblock \doi{10.1137/17M1122918}

\bibitem{davisUniversityFloridaSparse2011}
Davis, T.A., Hu, Y.: The {University of Florida} sparse matrix collection.
\newblock ACM Trans. Math. Softw. \textbf{38}(1) (2011).
\newblock \doi{10.1145/2049662.2049663}

\bibitem{edwardsKokkosEnablingManycore2014}
Edwards, H.C., Trott, C.R., Sunderland, D.: Kokkos: Enabling manycore
  performance portability through polymorphic memory access patterns.
\newblock J. Parallel Distr. Comput. \textbf{74}(12), 3202--3216 (2014).
\newblock \doi{10.1016/j.jpdc.2014.07.003}

\bibitem{giraudMixedprecisionPreconditionersParallel2008}
Giraud, L., Haidar, A., Watson, L.T.: Mixed-precision preconditioners in
  parallel domain decomposition solvers.
\newblock In: U.~Langer, M.~Discacciati, D.E. Keyes, O.B. Widlund, W.~Zulehner
  (eds.) Domain Decomposition Methods in Science and Engineering {{XVII}}, pp.
  357--364. {Springer Berlin Heidelberg}, {Berlin, Heidelberg} (2008).
\newblock \doi{10.1007/978-3-540-75199-1_44}

\bibitem{giraudLossOrthogonalityGramSchmidt2005}
Giraud, L., Langou, J., Rozloznik, M.: The loss of orthogonality in the
  {{Gram-Schmidt}} orthogonalization process.
\newblock Comput. Math. Appl. \textbf{50}(7), 1069--1075 (2005).
\newblock \doi{10.1016/j.camwa.2005.08.009}

\bibitem{grattonExploitingVariablePrecision2020}
Gratton, S., Simon, E., {Titley-Peloquin}, D., Toint, P.: Exploiting variable
  precision in {{GMRES}}.
\newblock SIAM J. Sci. Comput. (to appear)  (2020)

\bibitem{greenbaumNumericalBehaviourModified1997}
Greenbaum, A., Rozlo{\v z}n{\'i}k, M., Strako{\v s}, Z.: Numerical behaviour of
  the {Modified Gram-Schmidt} {{GMRES}} implementation.
\newblock Bit. Numer. Math. \textbf{37}(3), 706--719 (1997).
\newblock \doi{10.1007/BF02510248}

\bibitem{grutzmacherModularPrecisionFormat2018}
Gr{\"u}tzmacher, T., Anzt, H.: A modular precision format for decoupling
  arithmetic format and storage format.
\newblock In: Revised {{Selected Papers}}, pp. 434--443. {Turin, Italy} (2018).
\newblock \doi{10.1007/978-3-030-10549-5_34}

\bibitem{gustafsonBeatingFloatingPoint2017}
Gustafson, J.L., Yonemoto, I.T.: Beating floating point at its own game:
  {{Posit}} arithmetic.
\newblock Supercomput. Front. Innov. \textbf{4}(2), 71--86--86 (2017).
\newblock \doi{10.14529/jsfi170206}

\bibitem{liangErrorcontrolledLossyCompression2018}
Liang, X., Di, S., Tao, D., Li, S., Li, S., Guo, H., Chen, Z., Cappello, F.:
  Error-controlled lossy compression optimized for high compression ratios of
  scientific datasets.
\newblock In: 2018 {{IEEE}} International Conference on Big Data (Big Data),
  pp. 438--447. {IEEE} (2018).
\newblock \doi{10.1109/BigData.2018.8622520}

\bibitem{lindstromFixedrateCompressedFloatingpoint2014}
Lindstrom, P.: Fixed-rate compressed floating-point arrays.
\newblock IEEE Trans. Vis. Comput. Graph. \textbf{20}(12), 2674--2683 (2014).
\newblock \doi{10.1109/TVCG.2014.2346458}

\bibitem{paigeUsefulFormUnitary2009}
Paige, C.C.: A useful form of unitary matrix obtained from any sequence of unit
  2-norm n-vectors.
\newblock SIAM J. Matrix Anal. Appl. \textbf{31}(2), 565--583 (2009).
\newblock \doi{10.1137/080725167}

\bibitem{paigeModifiedGramSchmidtMGS2006}
Paige, C.C., Rozlozn{\'i}k, M., Strakos, Z.: Modified {{Gram-Schmidt}}
  ({{MGS}}), least squares, and backward stability of {{MGS-GMRES}}.
\newblock SIAM J. Matrix Anal. Appl. \textbf{28}(1), 264--284 (2006).
\newblock \doi{10.1137/050630416}

\bibitem{paigeResidualBackwardError2001}
Paige, C.C., Strakos, Z.: Residual and backward error bounds in minimum
  residual {Krylov} subspace methods.
\newblock SIAM J. Sci. Comput. \textbf{23}(6), 1898--1923 (2001).
\newblock \doi{10.1137/S1064827500381239}

\bibitem{saadIterativeMethodsSparse2003}
Saad, Y.: Iterative Methods for Sparse Linear Systems, second edn.
\newblock {SIAM Press}, {Philadelphia, PA, USA} (2003)

\bibitem{saadGMRESGeneralizedMinimal1986}
Saad, Y., Schultz, M.H.: {{GMRES}}: A generalized minimal residual algorithm
  for solving nonsymmetric linear systems.
\newblock SIAM J. Sci. Stat. Comput. \textbf{7}(3), 856--869 (1986).
\newblock \doi{10.1137/0907058}

\bibitem{vandervorstSuperlinearConvergenceBehaviour1993}
{Van der Vorst}, H.A., Vuik, C.: The superlinear convergence behaviour of
  {{GMRES}}.
\newblock J. Comput. Appl. Math. \textbf{48}(3), 327--341 (1993).
\newblock \doi{10.1016/0377-0427(93)90028-A}

\bibitem{wilkinsonRoundingErrorsAlgebraic1963}
Wilkinson, J.H.: Rounding Errors in Algebraic Processes.
\newblock {Prentice-Hall}, {Princeton, NJ, USA} (1963)

\end{thebibliography}
\end{document}